\documentclass[11pt]{article}
\usepackage{amsmath}
\usepackage{amssymb}
\usepackage{mathptmx}
\usepackage{amsfonts}
\usepackage[mathscr]{eucal}
\usepackage[dvips]{graphicx}
\usepackage{color}

\newcommand{\mbf}[1]{\ensuremath{\mathbf{#1}}}
\newcommand{\ms}[1]{\ensuremath{\mathscr{#1}}}

\newcommand{\tens}{\otimes}

\newcommand{\Min}{\textrm{min}}
\newcommand{\Max}{\textrm{max}}
\newcommand{\w}{\textrm{w}}

\newtheorem{theorem}{Theorem}

\newtheorem{corollary}[theorem]{Corollary}

\begin{document}

\title{A Hardy's  Uncertainty Principle Lemma     in  Weak  \\ Commutation   Relations  of Heisenberg-Lie Algebra}

\author{Toshimitsu TAKAESU }
\date{ }
\maketitle

\begin{center}
\textit{Faculty of Mathematics, Kyushu University,\\  Fukuoka, 812-8581, Japan }
\end{center}

\begin{quote}
\textbf{Abstract}.
  In this article  we consider linear operators satisfying   a generalized commutation relation of a type of the Heisenberg-Lie algebra. It is  proven that a generalized inequality  of the Hardy's uncertainty principle lemma follows. 
Its applications  to   time operators and abstract Dirac operators are  also investigated. \\
$\;$ 
$\qquad \qquad \; \; $ {\small time-energy uncertainty relation, Dirac operators, essential self-adjointness. } 
\end{quote}

\section{Introduction and Results}
In this article we investigate a norm-inequality of the linear operators
which obey  a  generalized weak commutation relation of  a type of the  Heisenberg-Lie algebra, and consider its application to   the theory of the time operator \cite{Mi01, Ar05}, and an abstract Dirac operator. 
Let   $\mbf{X} = \{X_j \}_{j=1}^N$, $\mbf{Y} = \{ Y_j \}_{j=1}^N$
 and $\mbf{Z} = \{Z_j \}_{j=1}^N$ be symmetric operators on a Hilbert space $\ms{H}$.  
The weak commutator of  operators $A$ and $B$ is defined   for $\psi  \,   \in \ms{D} (A) \cap \ms{D} (B)$ and 
 $\phi  \,   \in \ms{D} (A^\ast ) \cap \ms{D} (B^\ast) $ by 
\[
[ A  , B]^\w ( \phi  \, , \,  \psi  ) \; \; = \; \; ( A^\ast \phi  \,  , \, B  \psi   ) \; - \;  ( B^\ast \phi , \, A \psi  ). 
\]
Here the inner product has a linearity of $  (\eta , \alpha \psi + \beta\phi ) = \alpha (\eta , \psi )  + \beta ( \eta , \phi )$ for $\alpha , \beta \in \mbf{C}$. 
We assume that $(\mbf{X}, \mbf{Y}, \mbf{Z}) $ satisfies the following conditions. 
\begin{quote}
\textbf{(A.1)} $Z_j$,  $1 \leq j \leq  N $, is bounded operator. 
\end{quote}

\begin{quote}
\textbf{(A.2)} 
 Let 
  $\ms{D}_{\mbf{X}} = \cap_{j=1}^{N} \ms{D} (X_j )$ and $\ms{D}_{\mbf{Y}} = \cap_{j=1}^{N} \ms{D} (Y_j )$. 
It follows that  for $ \phi , \psi \in \ms{D}_{\mbf{X}} \cap \ms{D}_{\mbf{Y}}$, 
 \begin{align*}
 & [ X_j  , Y_l ]^\w ( \phi \, , \psi  )  \; = \; 
   \delta_{j, \, l } \, (\phi , \, i \, Z_j  \psi  )  , \\
 &[ X_j  , Z_l ]^\w ( \phi \, , \psi  )  \; = \;  [ Y_j  , Z_l ]^\w ( \phi \, , \psi  ) \; \; = \; \; 0 \\
 &  [ X_j  , X_l ]^\w ( \phi \, , \psi  ) \;  =  \; [ Y_j  , Y_l ]^\w ( \phi \, , \psi  ) \;  =\; [ Z_j  , Z_l ]^\w ( \phi \, , \psi  ) \; = \;  0 . 
 \end{align*}
 \end{quote}
Note that   $ \; [ Z_j  , Z_l ]  \psi \, =  \;  0 \;  $ follows for $\; \psi \in \ms{H}$, 
since $Z_j$, $j=1, \cdots, N$, is bounded.
 In this article we consider  an generalization of the inequality 
\begin{equation}
  \int_{\mbf{R}^N} \frac{1}{ |  \mbf{r}  |^2 }   | u  (\mbf{r}) |^2 d \mbf{r}  \; \; \leq 
\; \;   \frac{4}{(N-2)^2} \int_{\mbf{R}^N} \left|   \nabla u (\mbf{r}) \right|^2  d \mbf{r}  ,  \qquad \qquad  N\geq 3 .
  \notag
\end{equation}
This inequality is a basic one of Hardy's uncertainty principle inequalities. 
For Hardy's uncertainty inequalities, refer to   e.g.  \cite{ FS97, HJ94, Th04}.

$\; $ \\
Let us introduce  the additional conditions. 
\begin{quote}
\textbf{(A.3)}  $X_j$  is  self-adjoint for all $1 \leq j \leq  N $. 
\end{quote}
\begin{quote}
\textbf{(A.4)} $X_i$ and $Z_l$  strongly commutes for all $1 \leq j \leq  N $ and $1 \leq l \leq  N $ . 
\end{quote}

$\;$ \\ 
Sicne $Z_j  $, $j=1 , \cdots , N$, is bounded self-adjoint operator,   we can  set  $ \lambda_{\Min}  (\mbf{Z})  $ and 
 $ \lambda_{\Max}  (\mbf{Z}) $  by
\begin{align*}
& \lambda_{\Min}  (\mbf{Z}) \; \,  = \; \;      \min_{1\leq j \leq N }   \; \,  \inf  \sigma (Z_j)  ,  \\
& \lambda_{\Max}  (\mbf{Z})  \;  \, = \; \; \max_{1\leq j \leq N }  \; \,  \sup  \sigma (Z_j) , 
\end{align*}
where $\sigma (O ) $ denotes the spectrum of the operator $O$. \\

\begin{theorem}  \label{MaTheorem}
Assume \textbf{(A.1)}-\textbf{(A.4)}.  Let $ \Psi \in \ms{D} (|\mbf{X} |^{-1} ) \cap \ms{D}_{\mbf{X} }  \cap \ms{D}_{\mbf{Y}}$.  Then  the following \textbf{(1)} and \textbf{(2)} hold \\ 
 \textbf{(1)}  If 
 $ \;   N \lambda_{\Min} (\mbf{Z})  - 2 \lambda_{\Max}  (\mbf{Z})  > 0  $, it follows that 
\begin{equation}  
\left\| \frac{}{}  |\mbf{X}|^{-1} \Psi   \right\|^2    \; \leq \; \; 
\frac{4}{   \left( N \lambda_{\Min} (\mbf{Z})  - 2 \lambda_{\Max}  (\mbf{Z}) \right)^2} \, 
  \sum_{j=1}^N \left\|  \frac{}{}  Y_j  \, \Psi  \, \right\|^2  .    \label{G-Hardy1}
\end{equation}
\textbf{(2)}  If 
 $ \;   2 \lambda_{\Min} (\mbf{Z})  - N \lambda_{\Max}  (\mbf{Z})  > 0  $, it follows that 
\begin{equation}
\left\| \frac{}{}  |\mbf{X}|^{-1} \Psi   \right\|   \; \leq  \; \; 
\frac{4}{  \left(  2 \lambda_{\Min} (\mbf{Z})  - N \lambda_{\Max}  (\mbf{Z}) \right)^2 }
 \,  \sum_{j=1}^N  \, \left\|  \frac{}{}  Y_j  \, \Psi \, \right\|^2   .   \label{G-Hardy2}
\end{equation}
$\;$
\end{theorem}

$\; $ \\
Before proving Theorem \ref{MaTheorem}, let us consider the  replacement  of $\mbf{X}$ and $\mbf{Y}$ in Theorem 1.  
Let us introduce the following conditions substitute for \textbf{(A.3)} and \textbf{(A.4)}. 
\begin{quote}
\textbf{(A.5)}  $Y_j$  is  self-adjoint for all $1 \leq j \leq  N $. 
\end{quote}
\begin{quote}
\textbf{(A.6)} $Y_i$ and $Z_l$  strongly commutes for all $1 \leq j \leq  N $ and $1 \leq l \leq  N $ . 
\end{quote}
It is seen from  \textbf{(A.2)},  that 
\begin{equation}
\qquad \qquad \qquad 
[Y_j , \, X_l  ]^\w (\phi , \psi ) \; = \; \delta_{j, l} (\phi , \,  i \, (-Z_j ) \psi ) , \qquad \qquad \phi , \psi \in \ms{D}_{\mbf{X} }  \cap \ms{D}_{\mbf{Y}}  .
\end{equation}
Note  that $ \; \inf \sigma ( -  Z_j  ) = - \sup (Z_j  ) \; $  and $ \;  \sup (- Z_j) = - \inf \sigma (Z_j) \; $ follow. Then we obtain a following corollary : 

\begin{corollary}
 Assume \textbf{(A.1)}-\textbf{(A.2)} and \textbf{(A.5)}-\textbf{(A.6)}.  Let $ \Psi \in \ms{D} (|\mbf{Y} |^{-1} ) \cap \ms{D}_{\mbf{X} }  \cap \ms{D}_{\mbf{Y}}$.  Then  the following \textbf{(1)} and \textbf{(2)} hold. \\ 
\textbf{(1)}  If 
 $ \;   2 \lambda_{\Min} (\mbf{Z})  - N \lambda_{\Max}  (\mbf{Z})  > 0  $, it follows that 
\begin{equation}
\left\| \frac{}{}  |\mbf{Y}|^{-1} \Psi   \right\|   \; \leq  \; \; 
\frac{4}{  \left(  2 \lambda_{\Min} (\mbf{Z})  - N \lambda_{\Max}  (\mbf{Z}) \right)^2 }
 \,  \sum_{j=1}^N  \, \left\|  \frac{}{}  X_j  \, \Psi \, \right\|^2   .   \label{G-Hardy3}
\end{equation}
\textbf{(2)}  If 
 $ \;   N \lambda_{\Min} (\mbf{Z})  - 2 \lambda_{\Max}  (\mbf{Z})  > 0  $, it follows that 
\begin{equation}  
\left\| \frac{}{}  |\mbf{Y}|^{-1} \Psi   \right\|    \; \leq \; \; 
\frac{4}{\left(  N \lambda_{\Min} (\mbf{Z})  - 2 \lambda_{\Max}  (\mbf{Z})  \right)^2} \, 
  \sum_{j=1}^N \left\|  \frac{}{}  X_j  \, \Psi  \, \right\|^2  .    \label{G-Hardy4}
\end{equation}
$\;$ 
\end{corollary}

$\;$ \\
\textbf{{\large (Proof of Theorem \ref{MaTheorem})}} \\
 \textbf{(1)}Let $ \Psi \in \ms{D} (|\mbf{X} |^{-1} ) \cap \ms{D}_{\mbf{X} }  \cap \ms{D}_{\mbf{Y}} $. 
For  $\epsilon > 0 $ and $t>0 $,   it is seen that 
{\small
\begin{equation}
\left\| \left(  Y_j - it X_j (\mbf{X}^2 + \epsilon )^{-1}   \right) \Psi \right\|^2  \; \;   = \; \; \| Y_j \Psi  \|^2 \; \; 
  - it \, [  Y_j , \; X_j (\mbf{X} ^2  + \epsilon)^{-1} ]^\w ( \Psi , \, \Psi  ) \; \; + \; \; 
 t^2  \left\|  X_j (\mbf{X}^2 + \epsilon )^{-1}    \Psi \right\|^2 .  \label{a}
\end{equation}
}We see that
\begin{equation}
[  Y_j , \; X_j (\mbf{X}^2  + \epsilon )^{-1} ]^\w ( \Psi , \Psi  ) \; \; =
\; \; [  Y_j , \, X_j  ]^\w ( \Psi , \, (\mbf{X}^2  + \epsilon )^{-1} \Psi  ) \;  + 
 [  Y_j , \;  (\mbf{X}^2  + \epsilon )^{-1} ]^\w ( X_j \Psi , \Psi  )  \label{b} .
\end{equation}
From   \textbf{(A.2)} and  \textbf{(A.4)}, we obtain that
\begin{equation}
 [  Y_j , \, X_j  ]^\w ( \Psi , \, (\mbf{X}^2  + \epsilon )^{-1} \Psi  ) \; 
 =  \; -i \, ( (\mbf{X}^2  + \epsilon )^{-1/2} \Psi   , \, Z_j (\mbf{X}^2  + \epsilon )^{-1/2} \Psi )  . \label{c}
\end{equation}
Note that for  a symmetric operator $A$ and the non-negative symmetric operator $B$,  the resolvent formula \\ 
$
 [A, (B + \lambda )^{-1} ]^\w (v, u) =  [ B,  A]^\w ((B+ \lambda )^{-1} v, \, (B+ \lambda)^{-1} u ) $ for $ \lambda >0 $ follows.  Then by using this formura,  \textbf{(A.2)} and 
 \textbf{(A.4)} yield that
\begin{equation}
[  Y_j , \;  (\mbf{X}^2  + \epsilon )^{-1} ]^\w ( X_j \Psi , \Psi  ) \; = \; 
2i ( X_j ( \mbf{X}^2  + \epsilon )^{-1}  u \, , \, Z_j  X_j (\mbf{X}^2  + \epsilon )^{-1} u )
  \label{d}
\end{equation}
Since $\left\| \left(  Y_j - it X_j (\mbf{X}^2 + \epsilon )^{-1}   \right) u\right\|^2 \geq 0 \; $ and $ \; t>0$, we see from  (\ref{b}), (\ref{c}) and (\ref{d}) that 
\begin{align}
  &  \| Y_j \Psi  \|^2  \; \notag   \\
 &  
 \geq \;  -t^2 \left\|  X_j (\mbf{X}^2 + \epsilon )^{-1}   u\right\|^2 \; 
 +  t \, ( (\mbf{X}^2  + \epsilon )^{-1/2} \Psi  , \, Z_j (\mbf{X}^2  + \epsilon )^{-1/2} u ) - 2t ( X_j ( \mbf{X}^2  + \epsilon )^{-1}  u \, , \, Z_j  X_j (\mbf{X}^2  + \epsilon )^{-1} \Psi  )  \notag  \\
 & \geq \;  \left( \frac{}{} -t^2  - 2 t \lambda_{\Max}  (\mbf{Z})  \right)\left\|  X_j (\mbf{X}^2 + \epsilon )^{-1}   u\right\|^2 
 \; + t \lambda_{\Min}  (\mbf{Z}) \| (\mbf{X}^2  + \epsilon )^{-1/2} \Psi   \| 
   .  \label{e}
\end{align}
Then we have that 
\begin{equation}
  \sum_{j=1}^N\|   Y_j \, \Psi  \|^2  \;
 \geq \;  \left( \frac{}{} -t^2  - 2 t \lambda_{\Max}  (\mbf{Z})  \right)\left\|  |\mbf{X}| (\mbf{X}^2 + \epsilon )^{-1}   \Psi \right\|^2 
 \; + t N \lambda_{\Min}  (\mbf{Z}) \| (\mbf{X}^2  + \epsilon )^{-1/2} \Psi   \| 
   .  \label{f}
\end{equation}
Note that  $ \lim\limits_{\epsilon \to 0 }\left\|  |\mbf{X}| (\mbf{X}^2 + \epsilon )^{-1}    \Psi \right\|^2 
 \; = \;  \| |\mbf{X}|^{-1} \Psi \|   $ and $  \lim\limits_{\epsilon \to 0 }\| (\mbf{X}^2  + \epsilon )^{-1/2} \Psi  \|
\; = \; \| |\mbf{X}|^{-1} \Psi \| = 0  $ follow from  the spectral decomposition theorem. Then we have
 \begin{equation}
  \sum_{j=1}^N\|   Y_j \, \Psi  \|^2 \;
\; \geq  \; \;  \left(  -t^2    + (  N \lambda_{\Min} (\mbf{Z})  - 2 \lambda_{\Max}  (\mbf{Z}) )t \right)
\left\| \frac{}{}  |\mbf{X}|^{-1} \Psi    \right\|  . \label{gg}
 \end{equation}
 By taking  $ t =  \frac{  N \lambda_{\Min} (\mbf{Z})  - 2 \lambda_{\Max}  (\mbf{Z}) }{2} > 0  $ in the right side of  (\ref{gg}),  we obtain  \textbf{(1)}.  \\ 
 \textbf{(2)} By computing $\left\| \left(  Y_j  + it X_j (\mbf{X}^2 + \epsilon )^{-1}   \right) \Psi \right\|^2$ for $t>0 $ and $\epsilon >0$, in a similar way of \textbf{(1)}, we see that  
 \begin{align}
  &  \| Y_j \Psi  \|^2  \; \notag   \\
 &  
 \geq \;  -t^2 \left\|  X_j (\mbf{X}^2 + \epsilon )^{-1}   u\right\|^2 \; 
 -  t \, ( (\mbf{X}^2  + \epsilon )^{-1/2} \Psi  , \, Z_j (\mbf{X}^2  + \epsilon )^{-1/2} u ) + 2t ( X_j ( \mbf{X}^2  + \epsilon )^{-1}  u \, , \, Z_j  X_j (\mbf{X}^2  + \epsilon )^{-1} \Psi  )  \notag  \\
 & \geq \;  \left( \frac{}{} -t^2  + 2 t \lambda_{\Min}  (\mbf{Z})  \right)\left\|  X_j (\mbf{X}^2 + \epsilon )^{-1}   u\right\|^2 
 \; - t \lambda_{\Max}  (\mbf{Z}) \| (\mbf{X}^2  + \epsilon )^{-1/2} \Psi   \| 
   .  \label{e2}
\end{align}
 Then by taking $\epsilon \to 0 $   in the right side of  (\ref{e2}), it follows that  
  \begin{equation}
\sum_{j=1}^N\|   Y_j \, \Psi  \|^2    \;
\; \geq  \; \;  \left(  -t^2    +  (   \frac{}{} 2 \lambda_{\Min} (\mbf{Z})  - N \lambda_{\Max} (\mbf{Z}) )t \right)
\left\| \frac{}{}  |\mbf{X}|^{-1} \Psi    \right\|  . \label{g}
 \end{equation}
 By taking  $ t =  \frac{(  2 \lambda_{\Min} (\mbf{Z})  - N \lambda_{\Max}  (\mbf{Z}) )}{2} > 0  $ in (\ref{g}),   we obtain \textbf{(2)}. $\blacksquare$.

$\;$ \\  
\section{Applications}
 \subsection{Time-Energy Uncertainty inequality}
 In this subsection we consider an applicaion  to the theory of time operators \cite{Ar05, Mi01}. 
Let $H$, $T$, and $C$ be  linear operators on a Hilbert space $H$. It is said that 
  $H$ has the weak time operator $T$ with the uncommutative factor $C$ if $(H,T,C)$ satisfy the following conditions. 
  \begin{quote}
  \textbf{(T.1)} $H$ and $T$ are symmetric.
  \end{quote}
  \begin{quote}
  \textbf{(T.2)} C is bounded and self-adjoint. 
  \end{quote}
  \begin{quote}
  \textbf{(T.3)} It follows that for $\phi , \psi \in \ms{D} (H) \cap \ms{D} (T)$, 
  \[
  [ T , H]^\w (\phi , \psi  )\; = \;   \, ( \phi , C \psi ) .
  \]
  \end{quote}
  \begin{quote}
  \textbf{(T.4)}
  \[
  \delta_{C} \; := \inf_{\psi \in (\text{ker}C)^{\bot} \backslash \{ 0 \} }  
  \frac{\left|  (\Psi , C \,  \Psi ) \right| }{\| \psi \|^2} \; > \;  0 .
  \]
  \end{quote}
Assume that $(H,T,C) $ satisfies \textbf{(T.1)}-\textbf{(T.4)}.  
Then by using $ \| A u \| \, \|  B u \|  \geq | \text{Im} (Au, Bu ) | \geq \frac{1}{2}\left| \frac{}{}  [A, B]^\w (u, u) \right|$, 
it is seen that (H,T,C) satisfies the time-energy uncertainty inequality (\cite{Ar05}, Proposition4.1): 
\begin{equation} \qquad 
\frac{\left\| \frac{}{} \left( H - <H>_{\psi }  \right) \psi \right\|
\left\| \frac{}{} \left( T - <T>_{\psi }  \right) \psi \right\|}{\|\psi \|^2} \; \; \geq \; \; 
 \frac{\delta_{C}}{2}  ,   \qquad \psi \, \in \, \ms{D}(H) \cap \ms{D} (T) ,   \label{7/1.1}
\end{equation}
 where $<O>_{\psi} \; = \; ( \psi , O \psi )$. 
 From \textbf{(2)} in Theorem \ref{MaTheorem} and \textbf{(1)} in Corollary 2,  we obtain another type of the inequality between $T$ and  $H$ :

\begin{corollary} (\textbf{Time-Energy Uncertainty Inequalities})$\, $ \\ 
 Assume \textbf{(T.1)}-\textbf{(T.3)}. Then the following \textbf{(i)} and \textbf{( ii)} hold. \\
 \textbf{(i)} If $T$ is self-adjoint, $C$ and $T$ strongly commute, and $ \sup \sigma (C) < 2 \inf \sigma (C)$,
it follows that for   $ \psi \in \ms{D} (|T |^{-1} ) \cap \ms{D} (T)  \cap \ms{D} (H) $,
 \begin{equation}  
\left\| \frac{}{}  | T |^{-1} \psi   \right\|    \; \leq \; \; 
\frac{2}{  2 \inf \sigma (C) - \sup \sigma (C) } \, 
 \left\|  \frac{}{}  H    \, \Psi \right\|    .    \label{time-energy1}
\end{equation}
 \textbf{(ii)} If $H$ is self-adjoint, $C$ and $H$ strongly commute, and $ \sup \sigma (C) < 2 \inf \sigma (C)$,
it follows that for   $ \psi \in \ms{D} (|H |^{-1} ) \cap \ms{D} (H)  \cap \ms{D} (T) $,
 \begin{equation}  
\left\| \frac{}{}  | H |^{-1} \psi   \right\|    \; \leq \; \; 
\frac{2}{  2 \inf \sigma (C) - \sup \sigma (C) } \, 
 \left\|  \frac{}{}  T    \, \Psi \right\|    .    \label{time-energy1}
\end{equation}
 \end{corollary}

 \subsection{Abstract Dirac Operators with Coulomb Potential }
 Next tlt us consider the application to abstract Dirac operators. 
 We consider  the self-adjoint operators 
$ \mbf{P} \, = \, \{   P_j \}_{j=1}^N $ and $  \mbf{Q}  \, = \, \{ Q_j \}_{j=1}^N$ on a Hilbert space $\ms{H}$. Let us set      a  subspace $\ms{D} \subset   \cap_{j, l} \left( \ms{D} (P_j ) \cap \ms{D} (Q_l ) \right)$.  It is said that $(\ms{H}, \ms{D}, \mbf{P} , \mbf{Q} )_{N}$ is the weak representaion of the CCR with degree $N$, if  
 \ms{D} is dense in $\ms{H}$ and it follows that  for $\phi , \,  \psi \; \in  \ms{D} $,
\begin{align*}
& [  P_j  , \, Q_l   ]^\w  (\phi , \, \psi )\; = \;  i \delta_{j ,l }  (\phi , \psi ) ,   \\ 
& [  P_j  , \, P_l   ]^\w (\phi , \, \psi )\; = \;    [  Q_j  , \, Q_l   ]^\w (\phi , \, \psi )\; = \;    0  .
\end{align*}  
Let us define an  abstract Dirac operator as follows. 
Let $ ( \ms{H} , \ms{D} , \mbf{P} , \mbf{Q}  )_3$ be the weak representation of the CCR with degree three. 
      Let $\mbf{A} = \{ A_j \}_{j=1}^3$ and $B$ be the bounded self-adjoint operators on a Hilbert space $\ms{K}$.   Here $\mbf{A} = \{ A_j \}_{j=1}^3$ and $B$ satisfy the canonical anti-commutation relations $ \{ A_j , \, A_l \} = 2 \delta_{j, l } $, $  \{ A_j , B \} = 0 $,
 $B^2 = I_{\ms{K}} $ where $I_{\ms{K}}$ is the identity operator on $\ms{K}$.  The state Hilbert space space is defined by $ \ms{H}_{\textrm{Dirac}} \; = \; \ms{K} \tens \ms{H} $. 
The free abstract Dirac operator is defiend by 
\[
H_0 \; = \;  \sum_{j=1}^3  A_j \tens P_j  \; + \; B \tens  M  .
\]
Here we assume the following condition.
\begin{quote}
\textbf{(D.1)} $P_j$ and $P_l$  strongly commute for $1 \leq  j \leq 3 $, $ 1 \leq l \leq 3$. $\; P_j $,   $1 \leq j \leq 3$, and $M $ strongly commute.
\end{quote}
Then it is seen   that $H_0^2 \Psi = \,  \left( \frac{}{} \mbf{P}^2  + M^2 \right) \Psi $ for  $ \Psi \in \ms{D}$.
The abstract Dirac Operator 
  with the Coulomb potential  is defined by 
\[
  H (\kappa )\; = \; H_0 \; + \; \kappa I_{\ms{K}} \, \tens   |\mbf{Q}|^{-1} ,
  \]
 where $\kappa \in \mbf{R}$ is a parameter  called the  coupling constant.
We assume that the following condition
\begin{quote}
\textbf{(D.2)} It follows that $ \ms{D} 
\; \subset  \; \ms{D} ( \mbf{Q}|^{-1} ) $.
\end{quote}
Then  it follows   from \textbf{(1)} in Theorem \ref{MaTheorem} that  for
 $ \psi \in \ms{D} $,
 \[
 \|  I_{\ms{K}} \, \tens   |\mbf{Q}|^{-1} \psi \|^2  \;  \leq 
  \;   4 \sum_{j=1}^3 \| P_j \Psi \|^2  \; \leq 4 \, \| H_0 \Psi \|^2.  
 \]
  Hence by   the Kato-Rellich theorem, we obtaine the following corollary. \\
 
\begin{corollary}
Assume \textbf{(D.1)} and \textbf{(D.2)}. Then for $ | \kappa |  < \frac{1}{2}$, $H (\kappa )$ is essentially self-adjoint on $\ms{D}$.
 \end{corollary}
 $\;$ \\

$\quad$ \\
{\Large \textbf{Acknowledgments}} \\
It is pleasure to thank assistant professor Akito Suzuki and associate professor  Fumio Hiroshima for their advice and comments.
 

\begin{thebibliography}{99}
\small{
 \bibitem{Ar}
A.Arai, \textit{Mathematical principles of quantum phenamena}, Asakura-syoten, 2005.  (in japanese)
 \bibitem{Ar05}
A.Arai, Generalized weak Weyl relation and decay of quantum dynamics, \textit{Rev. Math. Phys.} \textbf{17} (2005)
 1071-1109.
\bibitem{Ar07}
A.Arai, Heisenberg Operators, invariant domains and Heisenberg equations of motion, \textit{Rev. Math. Phys}
 \textbf{19} (2007) 1045-1069.  
\bibitem{Ar07-2} 
 A.Arai, Spectrum of time operators. \textit{Lett. Math. Phys.} \textbf{80} 211-221 (2007).
\bibitem{FS97} 
G.B.Folland and A.Sitaram, The uncertainty principle : A mathematical survey, \textit{J. Fourier Anal. Appl.} \textbf{3} (1997) 207-238.
\bibitem{HJ94} 
V.Havin and B.Joricke, \textit{The uncertainty principle in harmonic analysis}, Springer 1994.
  \bibitem{Mi01}
M.Miyamoto, A generalized Weyl relation approach to the time operator and its connection to the survival probability,  \textit{J. Math. Phys} \textbf{42} (2001) 1038-1052. 
 \bibitem{MSME}
 J.G.Muga, R.S.Mayato,and I.L.Egsquiza (eds.), \textit{Time in quantum mechanics}. Springer 2002.
 \bibitem{PF95} P.Pfeifer and J. Fr\"olich, Generalized time-energy uncertainty relations and bounds on lifetimes of resonances, \textit{Rev. Mod. Phys.} \textbf{67} (1995) 759-779. 
 \bibitem{RS2}
M.Reed and B.Simon, \textit{Methods of Modern Mathematical Physics Vol.II},
Academic Press, 1979.
 \bibitem{Sc83} K.Schm\"{u}dgen, On the Heisenberg commutationn relationhn. I, J. Funct. Anal. \textbf{50} (1983) 8-49.
\bibitem{Tha}
 B. Thaller, \textit{The Dirac equation}, Springer, 1992.
\bibitem{Th04}
S.Thangavelu, \textit{An introduction to the uncertainty principle : Hardy's theorem on Lie groups}, Birkh\"auser, 2004.
}
\end{thebibliography}
\end{document}